# DISCUSSION: THE DANTZIG SELECTOR: STATISTICAL ESTIMATION WHEN $p$ IS MUCH LARGER THAN $n$[1]

By Ya'acov Ritov

*The Hebrew University of Jerusalem*

Candès and Tao, in an impressive and innovative paper, introduce an ingenious estimator. Their discussion brings back the standard $\ell_2$ loss function into the main focus of the "large $p$, small[er] $n$" discussion. We wish to present in this comment an apologia for using the prediction error criterion as the way to gauge the quality of the estimator in large-dimension models. This is not, however, a postmodernist essay on a cultural aspect of statistics. For this the reader may refer to the challenging discussion in Breiman ([1]). Our discussion is within the boundaries of the standard decision theory as applied to complex parameters.

The setup we consider is the standard structural point of view of regression (see Greenshtein and Ritov [2] for details). We observe an i.i.d. sample from the pair $(y, x)$, where $x$ is a $p$-dimensional random vector, while $y$ is real. We may, but do not need to, assume the linear structure $y = x'\beta_0 + z$, where the random variable $z$ is independent of $x$. The informal objective is to find a good estimator of $\beta$, so $z$ can be *defined* by being uncorrelated with the residuals $y - \beta_0'x$. At this stage of the discussion, an estimator cannot be said to be the best, since for that we should agree on an exact criterion, and unlike the situation with simple parametric models, an estimator will be asymptotically efficient when a specific risk function is considered, and not so if another criterion is applied.

The data we consider is $\mathcal{D}_n = \{(y_1, x_1), \dots, (y_n, x_n)\}$, a simple random sample from the distribution of $(y, x)$. We compute $\hat{\beta} = \hat{\beta}_n(\mathcal{D}_n)$. The prediction criterion compares $\hat{\beta}$ and $\beta_0$ not by a direct loss function, for example, $\mathcal{L}_2(\hat{\beta}, \beta_0) \equiv \|\hat{\beta} - \beta_0\|^2$, but indirectly by comparing the theoretical optimal $E(y - \beta_0'x)^2$ to the prediction performance of the estimator, $E((y - \hat{\beta}'x)^2|\mathcal{D}_n)$. The expectation is taken over the distribution of $(y, x)$ from which $\mathcal{D}_n$ is

Received January 2007.
[1]Supported in part by an ISF grant.

sampled. However, the prediction inefficiency is given by

$$PIE_2(\hat{\beta}) \equiv \mathrm{E}((y - \hat{\beta}'x)^2 | \mathcal{D}_n) - \mathrm{E}(y - \beta_0'x)^2$$

$$= (\hat{\beta} - \beta_0)'\Sigma_x(\hat{\beta} - \beta_0)$$

$$\equiv \|\hat{\beta} - \beta_0\|_{\Sigma_x}^2.$$

Now, let $\iota(\beta_0, \hat{\beta}) = \{i : \hat{\beta}_i^2 + \beta_{0i}^2 > 0\}$ and $\Sigma_x(\iota) = (\mathrm{cov}(x^i, x^j))_{i,j \in \iota}$, where $x = (x^1, \ldots, x^p)'$. With this notation

$$PIE(\hat{\beta}) = \|\hat{\beta} - \beta_0\|_{\Sigma_x(\iota)}^2,$$

and under the assumptions of Candès and Tao the $PIE_2$ and $\mathcal{L}_2$ criteria are comparable.

What are the $p$ regressors $x^1, \ldots, x^p$? There are two main possibilities. The first is they may be genuine $p$ explanatory variables, representing different measurements. Thus, in a particular investigation they can be, among other things, height, weight, income, socioeconomic status, gender, the number of visits to the supermarket, and so on. In the second extreme situation we start with very few explanatory variables (typically one), $u \in \mathbb{R}^k$, and the linear regression problem is defined in terms of $x^i = \psi_i(u)$, $i = 1, \ldots, p$. This is the situation we face in standard nonparametric regression techniques using wavelet techniques or cubic splines with fixed nodes, or in classification techniques like SVM (support vector machine) (where $x$ is defined explicitly using the "kernel trick").

When we are faced with the structural model with many different variables representing conceptually different properties (height and income, e.g.), it may make sense to assume that they are normal, but it is very unlikely that they are independent. Any strong assumption on the huge $n \times p$ matrix is hard to conceive. Just think about a $1000 \times 10,000$ matrix! Furthermore, the loss function $\mathcal{L}_2$ which makes sense in low dimensions, makes, by itself, little sense in high-dimensional spaces. In many cases it represents the average error of many different estimators of different quantities. The vector $\beta$ as a vector has very little meaning. It is just a collection of parameters. One may be interested in the impact of a single variable on the outcome (e.g., the number of previous visits to the supermarket), or of a small group of variables (e.g., representing the socioeconomic level), but no one has a simple interpretation for an eclectic list of 100 parameters. The only reason to consider the loss $p^{-1}\|\hat{\beta} - \beta_0\|^2$ is because it is the mean of the individual squared errors. However, there is no apparent reason why it is the arithmetic mean that should be taken. The prediction error is an "objective" way to find the right weighted mean of the individual errors, and it is strongly adapted to the particular situation at hand.

Suppose the relevant submatrix of the hat matrix, the one that is related to the active variables, has a high conditional number. The estimation problem is ill-posed, while the prediction problem is well-defined. If there is co-linearity, we may not be able to verify which one of the variables has a causal effect, and what is the proper value of each coefficient, but we can very well infer the aggregate impact of a group of variables. More problematic is the meaning of the proper value of a parameter. This can be defined as the population value. However, within the context of the "large $p$, small $n$ problem," the model is necessarily defined with respect to the sample size. The right model is the best that can be estimated with the given resources. With a larger sample size, we may want to use a completely different set of variables.

Consider now the other extreme in which the problem at hand is a nonparametric regression of $y$ on a univariate random variable $u$ whose distribution is unknown. The random vector $x$ is then $(\psi_1(u), \ldots, \psi_p(u))'$, where $\psi_1, \psi_2, \ldots$ is some basis of $L_2$. The assumption that the components of $x$ are normal seems now unreasonable. The assumption that the components are independent, or even uncorrelated, is very strong (at best, the $\psi_1, \psi_2, \ldots$ are orthonormal with respect to some a priori measure, e.g., Lebesgue, not the distribution of $u$). So, it is hard to see how much regularity can be assumed for the design matrix. Let $f_0(u) = \sum_{j=1}^p \beta_0^j \psi_j(u)$ and $\hat{f}(u) = \sum_{j=1}^p \hat{\beta}_0^j \psi_j(u)$. Then $\|\hat{\beta} - \beta_0\|^2 = \int (\hat{f}(u) - f_0(u))^2 \, du$ (assuming that the basis functions are Lebesgue orthogonal). This does make sense, but the prediction loss $\int (\hat{f}(u) - f_0(u))^2 \, dF_u(u)$ is still more reasonable.

However, restricting ourselves to orthonormal series (even with respect to the Lebesgue measure) may be too extreme. We want to invoke sparsity, which is essential to the large $p$ analysis, and sparsity depends on finding the proper $\psi_i$'s. The same function may be sparse in a given representation (e.g., the $\psi_i$'s are step functions or ramps) and not in other representations (e.g., when they are the Haar basis functions or step functions, resp.). Then $\|\hat{\beta} - \beta_0\|^2$ makes very little sense in terms of the estimated function $f$, but the prediction error is still exactly what one needs.

DEPARTMENT OF STATISTICS
THE HEBREW UNIVERSITY OF JERUSALEM
MOUNT SCOPUS, JERUSALEM 91905
ISRAEL
E-MAIL: yaacov.ritov@huji.ac.il